\documentclass{amsart}

\usepackage{amssymb}
\usepackage{amsmath}
\usepackage[all]{xy}

\newtheorem{theorem}[subsection]{Theorem}
\newtheorem*{ntheorem}{Theorem}
\newtheorem{cor}[subsection]{Corollary}

\newtheorem{lemma}[subsection]{Lemma}

\theoremstyle{remark}
\newtheorem{exa}[subsection]{Example}  
\newtheorem{rem}[subsection]{Remark}

\theoremstyle{definition}

\DeclareMathOperator{\BOX}{Box}
\DeclareMathOperator{\Hom}{Hom}

\DeclareMathOperator{\Span}{span}

\begin{document}

\title{Inequalities and Ehrhart $\delta$-Vectors}

\author{A. Stapledon}

\address{Department of Mathematics, University of Michigan, Ann Arbor, MI 48109, USA}
\email{astapldn@umich.edu}

\begin{abstract}
For any lattice polytope $P$, we consider an associated polynomial $\bar{\delta}_{P}(t)$ and describe its decomposition into 
a sum of two polynomials satisfying certain symmetry conditions. 
As a consequence,
we improve upon known inequalities  
satisfied by the coefficients of the Ehrhart $\delta$-vector of a lattice polytope. We also provide combinatorial proofs of two results of Stanley that were previously established using techniques from commutative algebra. Finally, we give a 
necessary numerical criterion for the existence of a regular unimodular lattice triangulation
of the boundary of a lattice polytope. 
\end{abstract}

\maketitle

\section{Introduction}

Let $N$ be a lattice of rank $n$ and set $N_{\mathbb{R}} = N \otimes_{\mathbb{Z}} \mathbb{R}$. 
Fix a $d$-dimensional lattice polytope $P$ in $N$ and, for each positive integer $m$, let $f_{P}(m)$ denote the 
number of lattice points in $mP$. It is a result of Ehrhart \cite{EhrLinearI, EhrLinearII} that
$f_{P}(m)$ is a polynomial in $m$ of degree $d$, called the 
\emph{Ehrhart polynomial} of $P$. 
The generating series of $f_{P}(m)$ can be written in the form
\begin{equation*}
\sum_{m \geq 0} f_{P}(m)t^{m} = \delta_{P}(t)/ (1 - t)^{d + 1},
\end{equation*}
where $\delta_{P}(t)$ is a polynomial of degree less than or equal to $d$,
called the \emph{(Ehrhart) $\delta$-polynomial} of $P$. 
If we write 
\begin{equation*}
\delta_{P}(t) = \delta_{0} + \delta_{1}t + \cdots +  \delta_{d}t^{d},
\end{equation*}
then $(\delta_{0}, \delta_{1}, \ldots, \delta_{d})$ is the \emph{(Ehrhart) $\delta$-vector} of $P$.
We will set 
$\delta_{i} = 0$ for $i < 0$ and $i > d$.  
It is a result of Stanley \cite{StaDecompositions} that the coefficients $\delta_{i}$ are non-negative integers. 
The degree $s$ of $\delta_{P}(t)$ is called the \emph{degree} of $P$
and $l = d + 1 - s$ is the \emph{codegree} of $P$. 
It is a consequence of Ehrhart Reciprocity that 
$l$ is the smallest 
positive integer such that $lP$ contains a lattice point in its relative interior (see, for example, \cite{HibAlgebraic}).
It is an open problem to characterise which vectors of non-negative integers
are $\delta$-vectors of a lattice polytope. Ideally, one would like a series of inequalities that are satisfied by exactly the $\delta$-vectors. 
We first summarise the current state of knowledge concerning inequalities and Ehrhart $\delta$-vectors. 

It follows from the definition that $\delta_{0} = 1$ and $\delta_{1} = f_{P}(1) - ( d + 1) = |P \cap N| - (d + 1)$. It is a consequence of Ehrhart Reciprocity that 
$\delta_{d}$ is the number of lattice points in the relative interior of $P$ (see, for example, \cite{HibAlgebraic}). Since $P$ has at least $d + 1$ vertices, we have the inequality $\delta_{1} \geq \delta_{d}$. 
We list the known inequalities satisfied by  the 
Ehrhart $\delta$-vector (c.f. \cite{BDDPSCoefficients}). 
\begin{equation}
\delta_{1} \geq \delta_{d}
\end{equation}
\begin{equation}\label{ineq1}
\delta_{0} + \delta_{1} 
+ \cdots  + \delta_{i + 1} \geq \delta_{d} + \delta_{d - 1} + \cdots + \delta_{d - i}
 \textrm{ for }     i = 0 , \ldots, \lfloor d/2     \rfloor - 1,
\end{equation}
\begin{equation}\label{ineq2}
\delta_{0} + \delta_{1} + \cdots + \delta_{i} \leq \delta_{s} + \delta_{s - 1} + \cdots + \delta_{s - i}
\textrm{ for }     i = 0 , \ldots, d,
\end{equation}
\begin{equation}\label{ineq3}
\textrm{ if } \delta_{d} \neq 0 \textrm{ then } 1 \leq \delta_{1} \leq \delta_{i} \textrm{ for }     i = 2 , \ldots, d - 1.
\end{equation}

Inequalities (\ref{ineq1}) and (\ref{ineq2}) were proved by Hibi in  \cite{HibSome} and Stanley in \cite{StaHilbert1} respectively. Both proofs are based  on commutative algebra. Inequality (\ref{ineq3}) was 
proved 
by Hibi in \cite{HibLower} using combinatorial methods. Recently, Henk and Tagami \cite{HTLower} produced examples showing that the analogue of (\ref{ineq3}) when $\delta_{d} = 0$ is false. 
That is, 
it is not true that $\delta_{1} \leq \delta_{i}$ for 
$i = 2, \ldots, s -1$. An explicit example is provided in Example \ref{cual} below. 

We improve upon these inequalities by proving the following result (Theorem \ref{yoke}). 
We remark that the proof is purely combinatorial. 
\begin{ntheorem}
Let $P$ be a $d$-dimensional lattice polytope of degree $s$ and codegree $l$.  The Ehrhart $\delta$-vector $(\delta_{0}, \ldots, \delta_{d})$
of $P$ satisfies the following inequalities.
\begin{equation}\label{I1}
\delta_{1} \geq \delta_{d},
\end{equation}
\begin{equation}\label{I2}
\delta_{2} + \cdots 
+ \delta_{i + 1} \geq \delta_{d - 1} + \cdots + \delta_{d - i}  \textrm{ for }     i = 0 , \ldots, \lfloor d/2     \rfloor - 1,
\end{equation}
\begin{equation}\label{I3}
\delta_{0} + \delta_{1} + \cdots + \delta_{i} \leq \delta_{s} + \delta_{s - 1} + \cdots + \delta_{s - i}
\textrm{ for }      i = 0, \ldots, d,
\end{equation}
\begin{equation}\label{I4}
\delta_{2 - l}  + \cdots + \delta_{0} +  \delta_{1} \leq \delta_{i} + \delta_{i - 1} + \cdots + \delta_{i - l + 1}  \textrm{ for }      
i = 2, \ldots, d - 1.
\end{equation}
\end{ntheorem}

\begin{rem}
Equality can be achieved in all the 
inequalities in the above theorem. For example, let $N$ be a lattice with basis $e_{1}, \ldots, e_{d}$ and 
let $P$ be the regular simplex with vertices $0,e_{1}, \ldots, e_{d}$. In this case, $\delta_{P}(t) = 1$
and each inequality above is an equality. 
\end{rem}

\begin{rem}
Observe that (\ref{I1}) and (\ref{I2}) imply that
\begin{equation}\label{ineq5}
\delta_{1} + \cdots + \delta_{i} + \delta_{i + 1} \geq \delta_{d} + \delta_{d - 1} + \cdots + \delta_{d - i},
\end{equation}
for $i = 0 , \ldots, \lfloor d/2     \rfloor - 1$. Since $\delta_{0} = 1$, we conclude that  (\ref{ineq1}) is always a 
strict inequality. We note that inequality (\ref{I2}) was suggested by Hibi in \cite{HibLower}.
\end{rem}

\begin{rem}
We can view the above result as providing, in particular, 
a combinatorial proof of Stanley's inequality (\ref{I3}).
\end{rem}

\begin{rem}
We claim that inequality (\ref{I4}) provides the correct generalisation of Hibi's inequality (\ref{ineq3}). Our contribution is 
to prove the cases when $l > 1$ and we refer the reader to \cite{HibLower} for a proof of  (\ref{ineq3}). 
In fact, 
in the proof of the above theorem
we show that (\ref{I4}) can be deduced from (\ref{ineq3}), (\ref{I2}) and (\ref{I3}). 
\end{rem}

In order to prove this result, we 
 consider the polynomial
\[
\bar{\delta}_{P}(t) = (1 + t + \cdots + t^{l - 1})\delta_{P}(t)
\]
and use a result of Payne (Theorem 2 \cite{PayEhrhart}) to establish the following decomposition theorem (Theorem \ref{goblue}).
This generalises a result of Betke and McMullen for lattice polytopes containing a lattice point in their relative interior 
(Theorem 5 \cite{BMLattice}). 

\begin{ntheorem}
The polynomial $\bar{\delta}_{P}(t)$ has a unique decomposition 
\begin{equation*}
\bar{\delta}_{P}(t) = a(t) + t^{l}b(t),
\end{equation*}
where $a(t)$ and $b(t)$  are polynomials with integer coefficients satisfying
$a(t) = t^{d}a(t^{-1})$ and $b(t) = t^{d - l}b(t^{-1})$.
Moreover, the coefficients of $b(t)$ are non-negative and, if $a_{i}$ denotes the coefficient of $t^{i}$ in $a(t)$, then
\begin{equation*}
1 = a_{0} \leq a_{1} \leq a_{i},
\end{equation*}
for $i = 2, \ldots, d - 1$.
\end{ntheorem}

Using some elementary arguments we show that our desired inequalities are equivalent to certain conditions on the coefficients
of $\bar{\delta}_{P}(t)$, $a(t)$ and $b(t)$ (Lemma \ref{magnify}) and hence are a consequence of the above theorem. 

We also consider the following result of Stanley (Theorem 4.4 \cite{StaHilbert2}), which was proved using commutative algebra. 
\begin{ntheorem}
If $P$ is a lattice polytope of degree $s$ and codegree $l$, then $\delta_{P}(t) = t^{s} \delta_{P}(t^{-1})$ if and only if 
$lP$ is a translate of a reflexive polytope.  
\end{ntheorem}
We show that this result is a consequence of the above decomposition of $\bar{\delta}_{P}(t)$, thus providing a combinatorial
proof of Stanley's theorem (Corollary \ref{tata}). A different combinatorial proof in the case $l = 1$ is provided in \cite{HibDual}.

Recent work of Athanasiadis \cite{AthVectors, AthEhrhart} relates the existence of certain triangulations of a lattice polytope $P$ to inequalities satisfied by its $\delta$-vector.  
\begin{ntheorem}[Theorem 1.3 \cite{AthVectors}]
Let $P$ be a $d$-dimensional lattice polytope. If $P$ admits a regular unimodular lattice triangulation, then 
\begin{equation}\label{un}
\delta_{i + 1} \geq \delta_{d - i} \textrm{ for } 
i = 0, \ldots, \lfloor d/2 \rfloor - 1,
\end{equation}
\begin{equation}
\delta_{\lfloor (d + 1)/2 \rfloor} \geq \cdots \geq \delta_{d - 1} \geq \delta_{d},
\end{equation}
\begin{equation}\label{tres}
\delta_{i} \leq \binom{\delta_{1} + i - 1}{i}
\textrm{ for } i = 0, \ldots, d. \end{equation}
\end{ntheorem}

As a corollary of our decomposition of $\bar{\delta}_{P}(t)$, we deduce the following theorem (Theorem \ref{dwarf}).



\begin{ntheorem}
Let $P$ be a $d$-dimensional lattice polytope. If the boundary of $P$ 
admits a regular unimodular lattice triangulation, then 
\begin{equation}\label{cuatro}
\delta_{i + 1} \geq \delta_{d - i} \end{equation} 
\begin{equation}\label{cinco} 
\delta_{0} + \cdots + \delta_{i + 1} \leq \delta_{d} + \cdots + \delta_{d - i} + \binom{\delta_{1} - \delta_{d} + i + 1}{i + 1}, \end{equation}
for 
$i = 0, \ldots, \lfloor d/2 \rfloor - 1$.
\end{ntheorem}

We note that (\ref{cuatro}) provides a generalisation of (\ref{un}) and that (\ref{cinco}) may be viewed as an analogue of (\ref{tres}).   
We remark that the method of proof is different to that of Athanasiadis.

\medskip

I would like to thank my advisor Mircea Musta\c t\v a for all his help.   I would also like to thank Sam Payne for some valuable feedback. 
The author was supported by
Mircea Musta\c t\v a's Packard Fellowship and
by an Eleanor Sophia Wood
travelling scholarship from the University of Sydney.

\section{Inequalities and Ehrhart $\delta$-Vectors}



We will use the definitions and notation from the introduction throughout the paper. 
Our main object of study will be the polynomial
\[
\bar{\delta}_{P}(t) = (1 + t + \cdots + t^{l - 1})\delta_{P}(t). 
\]
Since $\delta_{P}(t)$ has degree $s$ and non-negative integer coefficients, it follows that $\bar{\delta}_{P}(t)$ 
has degree $d$ and non-negative integer coefficients. In fact, we will show that $\bar{\delta}_{P}(t)$ has positive integer
coefficients (Theorem \ref{goblue}). 
Observe that we can recover $\delta_{P}(t)$ from $\bar{\delta}_{P}(t)$ if we know the 
codegree $l$ of $P$.
If we write 
\[
\bar{\delta}_{P}(t) = \bar{\delta}_{0} + \bar{\delta}_{1}t + \cdots +  \bar{\delta}_{d}t^{d}, 
\]
then 
\begin{equation}\label{pinochet}
\bar{\delta}_{i} = \delta_{i} + \delta_{i - 1} + \cdots + \delta_{i - l + 1},
\end{equation}
for $i = 0, \ldots, d$. Note that $\bar{\delta}_{0} = 1$ and $\bar{\delta}_{d} = \delta_{s}$. 

\begin{exa}
Let $N$ be a lattice with basis $e_{1}, \ldots, e_{d}$ and let $P$ be the \emph{standard simplex} with vertices 
$0, e_{1}, \ldots, e_{d}$. It can be shown that $\delta_{P}(t) = 1$ and hence $\bar{\delta}_{P}(t) = 1 + t + \cdots + t^{d}$.
On the other hand, if $Q$ is the \emph{standard reflexive simplex} with vertices 
$e_{1}, \ldots, e_{d}$ and $-e_{1} - \cdots - e_{d}$, then $\bar{\delta}_{Q}(t) = \delta_{Q}(t) =  1 + t + \cdots + t^{d}$. We conclude that
$\bar{\delta}_{P}(t)$ does not determine $\delta_{P}(t)$.
\end{exa}

\begin{rem}\label{tiger}
We can interpret $\bar{\delta}_{P}(t)$ as the Ehrhart $\delta$-vector of a $(d + l)$-dimensional polytope. More specifically, 
let $Q$ be the standard reflexive simplex of dimension $l - 1$ in a lattice $M$ as above. 
Henk and Tagami \cite{HTLower} defined $P \otimes Q$ to be the convex hull in $(N \times M \times \mathbb{Z})_{\mathbb{R}}$ of
$P \times \{0 \} \times \{ 0 \}$ and $\{ 0 \} \times Q \times \{ 1 \}$. By Lemma 1.3 in \cite{HTLower}, $P \otimes Q$ is a $(d + l)$-dimensional lattice
polytope with Ehrhart $\delta$-vector
$\delta_{P \otimes Q}(t) = \delta_{P}(t)\delta_{Q}(t) = \delta_{P}(t)( 1 + t + \cdots + t^{l - 1}) = \bar{\delta}_{P}(t)$
\end{rem}
Our main objects of study
will be the polynomials $a(t)$ and $b(t)$ in the following elementary lemma.

\begin{lemma}\label{badgers}
The polynomial $\bar{\delta}_{P}(t)$ has a unique decomposition 
\begin{equation}\label{pog}
\bar{\delta}_{P}(t) = a(t) + t^{l}b(t),
\end{equation}
where $a(t)$ and $b(t)$  are polynomials with integer coefficients satisfying
$a(t) = t^{d}a(t^{-1})$ and $b(t) = t^{d - l}b(t^{-1})$.
\end{lemma}
\begin{proof}
Let $a_{i}$ and $b_{i}$ denote the coefficients of $t^{i}$ in $a(t)$ and $b(t)$ respectively, and set
\begin{equation}\label{coeff1}
a_{i + 1} = \delta_{0} + \cdots + \delta_{i + 1} - \delta_{d} - \cdots - \delta_{d - i}, 
\end{equation}
\begin{equation}\label{coeff2}
b_{i} = -\delta_{0} - \cdots - \delta_{i} + \delta_{s} + \cdots + \delta_{s - i}.
\end{equation}
We compute, using (\ref{pinochet}) and since 
$s + l = d + 1$, 
\begin{align*}
a_{i} + b_{i - l} &= \delta_{0} + \cdots + \delta_{i} - \delta_{d} - \cdots - \delta_{d - i + 1} 
-\delta_{0} - \cdots - \delta_{i - l} + \delta_{s} + \cdots + \delta_{s - i+ l} \\
&= \delta_{i - l + 1} + \cdots + \delta_{i} = \bar{\delta}_{i}, \\
a_{i} - a_{d - i} &= \delta_{0} + \cdots + \delta_{i} - \delta_{d} - \cdots - \delta_{d - i + 1}
- \delta_{0} -  \cdots - \delta_{d - i} + \delta_{d} + \cdots + \delta_{i + 1} \\
&= 0, \\
b_{i} - b_{d - l - i} &= -\delta_{0} - \cdots - \delta_{i} + \delta_{s} + \cdots + \delta_{s - i}
+ \delta_{0} + \cdots + \delta_{s - i - 1} - \delta_{s} - \cdots - \delta_{i + 1} \\
&= 0,
\end{align*}
for $i = 0, \ldots, d$. Hence we obtain our desired decomposition and one easily verifies the uniqueness assertion. 
\end{proof}


\begin{exa}\label{cual}
Let $N$ be a lattice with basis $e_{1}, \ldots, e_{5}$ and let $P$ be the $5$-dimensional lattice polytope
with vertices $0$, $e_{1}$, $e_{1} + e_{2}$,
$e_{2} + 2e_{3}$, $3e_{4} + e_{5}$ and $e_{5}$. Henk and Tagami showed that 
$\delta_{P}(t) = (1 + t^{2})(1 + 2t) = 1 + 2t + t^{2} + 2t^{3}$ (Example 1.1 in \cite{HTLower}). It follows that 
$s = l = 3$ and $\bar{\delta}_{P}(t) = 1 + 3t + 4t^{2} + 5t^{3} + 3t^{4} + 2t^{5}$. We calculate that
$a(t) = 1 + 3t + 4t^{2} + 4t^{3} + 3t^{4} + t^{5}$ and $b(t) = 1 + 0t + t^{2}$. 
\end{exa}


We may view our proposed inequalities on the coefficients of the Ehrhart $\delta$-vector as conditions on the 
coefficients of $\bar{\delta}_{P}(t)$, $a(t)$ and $b(t)$.

\begin{lemma}\label{magnify}
With the notation above,
\begin{enumerate}
\item[-] Inequality (\ref{ineq1}) holds if and only if the coefficients of $a(t)$ are non-negative.
\item[-] Inequality (\ref{ineq2}) holds if and only if the coefficients of $b(t)$ are non-negative. 
\item[-] Inequality (\ref{I2}) holds if and only if $a_{1} \leq a_{i}$ for $i = 2, \ldots, d - 1$.
\item[-] Inequality (\ref{I3}) holds if and only if the coefficients of $b(t)$ are non-negative.
\item[-] Inequality (\ref{I4}) holds if and only if $\bar{\delta}_{1} \leq \bar{\delta}_{i}$ for $i = 2, \ldots, d - 1$.
\item[-] Inequality (\ref{ineq5}) holds if and only if the coefficients of $a(t)$ are positive.
\end{enumerate}
\end{lemma}
\begin{proof}
The result follows by substituting (\ref{pinochet}),(\ref{coeff1}) and (\ref{coeff2})
into the right hand sides of the above statements. 
\end{proof}

\begin{rem}\label{barnes}
The coefficients of $a(t)$ are \emph{unimodal} if $a_{0} = 1 \leq a_{1} \leq \cdots \leq a_{\lfloor d/2 \rfloor}$. It 
follows from (\ref{coeff1}) that $a_{i + 1} - a_{i} = \delta_{i + 1} - \delta_{d - i}$ for all $i$. 
Hence the coefficients of $a(t)$
are unimodal if and only if $\delta_{i + 1} \geq \delta_{d - i}$ for $i = 0, \ldots,  \lfloor d/2 \rfloor - 1$. 
In Remark \ref{bonjovi}, we show that the coefficients of $a(t)$ are unimodal for $d \leq 5$. 
\end{rem}

\begin{rem}\label{smh}
A lattice polytope $P$ is \emph{reflexive} if 
the origin is the unique  lattice point in its relative interior and each facet $F$ of $P$ 
has the form $F  = \{ v \in P \mid \langle u , v \rangle = -1 \}$,
for some $u \in \Hom (N, \mathbb{Z})$. Equivalently, $P$ is reflexive if 
it contains the origin in its 
relative interior
and, for every positive integer $m$, every non-zero lattice point in $mP$ lies on $\partial(n P)$ for a unique positive integer 
$n \leq m$. 
It is a result of Hibi \cite{HibDual} that $\delta_{P}(t) = t^{d}\delta_{P}(t^{-1})$ if and only if $P$ is a translate of a reflexive polytope (c.f. Corollary \ref{tata}).  
We see from Lemma \ref{badgers} that $\delta_{P}(t) = t^{d}\delta_{P}(t^{-1})$ if and only if $\delta_{P}(t)= \bar{\delta}_{P}(t) = a(t)$.  
Payne and Musta\c t\v a gave examples of reflexive polytopes where the coefficients of $\delta_{P}(t) = a(t)$ are not unimodal 
\cite{MPEhrhart}.
Further examples are given by Payne for all $d \geq 6$ \cite{PayEhrhart}. 
\end{rem}

\begin{rem}
It follows from (\ref{coeff2}) that $b_{i + 1} - b_{i} = \delta_{s - (i + 1)} - \delta_{i + 1}$ for all $i$. 
Hence the coefficients of $b(t)$ are unimodal if and only if $\delta_{i} \leq \delta_{s - i}$ for $i = 1, \ldots, \lfloor (s - 1)/2
\rfloor$. We see from Example \ref{cual} that the coefficients of $b(t)$ are not necessarily unimodal. 
\end{rem}

Our next goal is to express $\bar{\delta}_{P}(t)$ as a sum of shifted $h$-vectors, using a result of Payne 
(Theorem 1.2 \cite{PayEhrhart}). 
We first fix a lattice triangulation $\mathcal{T}$ of $\partial P$ and recall 
what it means for $\mathcal{T}$ to be regular. 
Translate $P$ by an element of $N_{\mathbb{Q}}$ so that the origin lies in its interior
and let $\Sigma$ denote the fan over the faces of $\mathcal{T}$. Then $\mathcal{T}$ is \emph{regular} if $\Sigma$ can be realised as the fan over the faces of a  rational polytope.
Equivalently, $\mathcal{T}$ is regular if the toric variety $X(\Sigma)$ is projective.
We may always choose $\mathcal{T}$ to be a regular triangulation (see, for example, \cite{AthVectors}).
We regard the empty face as a face of $\mathcal{T}$ of dimension $-1$.
For each face $F$ of $\mathcal{T}$, consider the $h$-vector of $F$,
\begin{equation*} \label{change}
h_{F}(t) = \sum_{F \subseteq G} t^{\dim G - \dim F} (1 - t)^{d - 1 - \dim G}.
\end{equation*}
We recall the following well-known lemma and outline a geometric proof. We refer the reader to \cite{FulIntroduction} for more 
details (see also \cite{StaNumber}). 

\begin{lemma}\label{patriots}
Let $\mathcal{T}$ be a regular lattice triangulation of $\partial P$. 
If $F$ is a face of $\mathcal{T}$, then the $h$-vector of $F$ 
is a polynomial of degree $d - 1 - \dim F$ with symmetric, unimodal integer coefficients.
That is, if $h_{i}$ denotes the coefficient of $t^{i}$ in $h_{F}(t)$, 
then
$h_{i} = h_{d - 1 - \dim F - i}$ for all $i$ and $1 = h_{0} \leq h_{1} \leq \cdots \leq h_{\lfloor (d - 1 - \dim F)/2 \rfloor}$. The coefficients satisfy the Upper Bound Theorem,
\[
h_{i} \leq \binom{h_{1} + i - 1}{i},
\]
for $i = 1, \ldots, d - 1 - \dim F$. 
\end{lemma}
\begin{proof}
As above, translate $P$ by an element of $N_{\mathbb{Q}}$ so that the origin lies in its interior
and let $\Sigma$ denote the fan over the faces of $\mathcal{T}$. 
For each face $F$ of $\tau$, let $\Span F$ denote the smallest
linear subspace of $N_{\mathbb{R}}$ containing $F$ and let $\Sigma_{F}$ be the complete fan in $N_{\mathbb{R}}/ \Span F$ whose cones 
are 
the projections of the cones in $\Sigma$ containing $F$.
We can interpret $h_{i}$ as the dimension of the $2i^{\textrm{th}}$ cohomology group of the projective toric variety $X(\Sigma_{F})$. 
The symmetry of the $h_{i}$ follows from Poincar\'e
Duality on   $X(\Sigma_{F})$, while unimodality follows from the Hard Lefschetz Theorem. The cohomology ring $H^{*}( X(\Sigma_{F}), \mathbb{Q} )$ is isomorphic to 
the quotient of a polynomial ring in 
$h_{1}$ variables of degree $2$ and hence $h_{i}$ is bounded by the number of monomials of degree $i$
 in $h_{1}$ variables of degree $1$. 
\end{proof}

Recall that $l$ is the smallest positive integer such that $lP$ contains a lattice point in its relative interior
and fix a lattice point $\bar{v}$ in $lP \smallsetminus \partial(lP)$. 
Let $N' = N \times \mathbb{Z}$ and let $u: N' \rightarrow \mathbb{Z}$ denote the projection onto the second factor. 
We write $\sigma$ for the cone over $P \times \{ 1 \}$ in $N'_{\mathbb{R}}$ and $\rho$ for the ray through $(\bar{v}, l)$. 
For each face $F$ of $\mathcal{T}$, 
let $\sigma_{F}$ denote the cone over $F$ and 
let $\sigma_{F}'$ denote the cone generated by $\sigma_{F}$ and $\rho$. 
The empty face corresponds to the origin and $\rho$ respectively. 
The union of such cones forms a simplicial fan $\triangle$
refining $\sigma$. 
For each non-zero cone $\tau$ in $\triangle$, with primitive integer generators $v_{1}, \ldots, v_{r}$, we consider the open parallelopiped
\[
\BOX(\tau) = \{  a_{1}v_{1} + \cdots + a_{r}v_{r} \mid 0 < a_{i} < 1    \},
\]
and observe that we have an involution
\[
\iota: \BOX(\tau) \cap N' \rightarrow \BOX(\tau) \cap N'
\]
\[
\iota(a_{1}v_{1} + \cdots + a_{r}v_{r}) = (1 - a_{1})v_{1} + \cdots + (1 - a_{r})v_{r}.
\]
We also set $\BOX(\{ 0 \}) = \{ 0 \}$ and $\iota(0) = 0$. Observe that $\BOX(\rho) \cap N' = \emptyset$. 
For each face $F$ of $\mathcal{T}$, we define 
\[ B_{F}(t) = 
\sum_{ v \in \BOX(\sigma_{F})\cap N'} t^{u(v)}
\] 
\[ B_{F}'(t) = 
\sum_{ v \in \BOX(\sigma_{F}')\cap N'} t^{u(v)}. 
\]
If $\BOX(\sigma_{F}) \cap N' = \emptyset$ or $\BOX(\sigma_{F}') \cap N' = \emptyset$ then we define 
$B_{F}(t) = 0$ or $B_{F}'(t) = 0$ respectively. 
For example, when $F$ is the empty face, $B_{F}(t) = 1$ and 
$B_{F}'(t) = 0$. 
 We will need the following lemma.

\begin{lemma}\label{ball}
For each face $F$ of $\mathcal{T}$, $B_{F}(t) = t^{\dim F + 1}B_{F}(t^{-1})$ and 
$B_{F}'(t) = t^{\dim F + l + 1}B_{F}'(t^{-1})$.
\end{lemma}
\begin{proof}
Using the involution $\iota$ above, 
\[
t^{\dim F + 1}B_{F}(t^{-1}) = \sum_{ v \in \BOX(\sigma_{F})\cap N'} t^{\dim F + 1 - u(v)} = 
\sum_{ v \in \BOX(\sigma_{F})\cap N'} t^{u(\iota(v))} = B_{F}(t). 
\]
Similarly, 
\[
t^{\dim F + l + 1}B_{F}'(t^{-1}) = \sum_{ v \in \BOX(\sigma_{F}')\cap N'} t^{\dim F + 1 + l - u(v)}
= \sum_{ v \in \BOX(\sigma_{F}')\cap N'} t^{u(\iota(v))} = B_{F}'(t).\]
\end{proof}

Consider any element $v$ in $\sigma \cap N'$ and let $G$ be the smallest face of $\mathcal{T}$ such that 
$v$ lies in $\sigma_{G}'$. Set $r = \dim G + 1$ and let $v_{1}, \ldots, v_{r}$ denote the vertices of $G$. Then $v$ can be uniquely written in the form
\begin{equation}\label{Bombay}
v = \{ v \} + \sum_{ (v_{i}, 1) \notin \tau} (v_{i},1) + w,
\end{equation}
where $\{ v \}$ lies in $\BOX(\tau) \cap N'$, for some subcone
$\tau$ of $\sigma_{G}'$, 
and $w$ is a  
non-negative 
integer sum of $(v_{1}, 1), \ldots, (v_{r}, 1)$ and $(\bar{v}, l)$. 
If we write $w = \sum_{i = 1}^{r} a_{i}(v_{i},1) + a_{r + 1}(\bar{v},l)$, 
for some non-negative integers $a_{1}, \ldots, a_{r + 1}$, then
\begin{equation}\label{ams}
u(v) = u(\{ v \}) + \dim G - \dim F + \sum_{i = 1}^{r} a_{i} + a_{r + 1}l. 
\end{equation}
Conversely, given 
$\tilde{v}$ in $\BOX(\tau) \cap N'$, for some $\tau \subseteq  \sigma_{G}'$, and $w$ a  
non-negative 
integer sum of $(v_{1}, 1), \ldots, (v_{r}, 1)$ and $(\bar{v}, l)$, then 
$v =  \tilde{v} + \sum_{ (v_{i}, 1) \notin \tau} (v_{i},1) + w$ lies in $\sigma \cap N'$ and $G$ is the smallest face of 
$\mathcal{T}$ such that 
$v$ lies in $\sigma_{G}'$.


\begin{rem}\label{quo}
With the above notation, observe that $\tau = \sigma_{F}$ for some $F \subseteq G$ if and only if $v$ lies on a translate of 
$\partial \sigma$ by a non-negative multiple of $(\bar{v}, l)$. Note that if $\tau = \sigma_{F}'$ for some 
(necessarily non-empty) $F \subseteq G$, 
then $\{v + nv_{1} \} = \{ v \}$ and $u(v + nv_{1}) = u(v) + n$, for any non-negative integer $n$. We conclude that
$B_{F}'(t) = 0$ for all faces $F$ of $\mathcal{T}$ if and only if every element $v$ in $\sigma \cap (N \times l\mathbb{Z})$
can be written as the sum of an element of $\partial(mlP) \times \{ ml \}$ and $m'(\bar{v}, l)$,
for some non-negative integers $m$ and $m'$. 
\end{rem}

The generating series of $f_{P}(m)$ can be written as $\sum_{v \in \sigma \cap N'} t^{u(v)}$. 
Payne described this sum by considering the contributions of all $v$ in $\sigma \cap N'$ with a 
fixed $\{ v \} \in \BOX(\tau)$.  
We have the following application of Theorem 1.2 in \cite{PayEhrhart}. We recall the proof in this situation for the convenience 
of the reader.

\begin{lemma}\label{upper}
With the notation above,
\[
\bar{\delta}_{P}(t) = \sum_{F \in \mathcal{T}} (B_{F}(t) + B_{F}'(t)) h_{F}(t).
\]
\end{lemma}
\begin{proof}
Using (\ref{Bombay}) and (\ref{ams}), 
we compute
\begin{align*}
\bar{\delta}_{P}(t) &= (1 + t + \cdots + t^{l - 1}) \delta_{P}(t) = (1 - t^{l})(1 - t)^{d} \sum_{v \in \sigma \cap N'} t^{u(v)} \\
&= (1 - t)^{d} \sum_{F \in \mathcal{T}} (B_{F}(t) + B_{F}'(t) ) \sum_{ F \subseteq G} t^{\dim G - \dim F}/(1 - t)^{\dim G + 1}   \\
&=  \sum_{F \in \mathcal{T}} (B_{F}(t) + B_{F}'(t)) h_{F}(t).
\end{align*}

\end{proof}

\begin{rem}\label{fruit}
We can write $\bar{\delta}_{P}(t) = (1 - t)^{d + 1} \sum_{v \in \sigma \cap N'} (1 + t + \cdots + t^{l - 1})t^{u(v)}$.
Ehrhart Reciprocity states that, for any positive integer $m$, $f_{P}(-m)$ is $(-1)^{d}$ times the number of lattice points in the 
relative interior of $mP$ 
(see, for example, \cite{HibAlgebraic}). 
Hence $f_{P}(-1) = \cdots = f_{P}(1 - l) = 0$ and  
the generating series of the polynomial $f_{P}(m) + f_{P}(m - 1) + \cdots + f_{P}(m - l + 1)$ has the form
$\bar{\delta}_{P}(t)/(1 - t)^{d + 1}$. 
\end{rem}

We will now prove our first main result. When $s = d$, $\bar{\delta}_{P}(t) = \delta_{P}(t)$ and the theorem below is due to Betke and McMullen (Theorem 5 \cite{BMLattice}). This case was also proved in Remark 3.5 in \cite{YoWeightI}.

\begin{theorem}\label{goblue}
The polynomial $\bar{\delta}_{P}(t)$ has a unique decomposition 
\begin{equation*}
\bar{\delta}_{P}(t) = a(t) + t^{l}b(t),
\end{equation*}
where $a(t)$ and $b(t)$  are polynomials with integer coefficients satisfying
$a(t) = t^{d}a(t^{-1})$ and $b(t) = t^{d - l}b(t^{-1})$.
Moreover, the coefficients of $b(t)$ are non-negative and, if $a_{i}$ denotes the coefficient of $t^{i}$ in $a(t)$, then
\begin{equation}\label{lanquin}
1 = a_{0} \leq a_{1} \leq a_{i},
\end{equation}
for $i = 2, \ldots, d - 1$.
\end{theorem}

\begin{proof}
Let $\mathcal{T}$ be a regular lattice triangulation of $\partial P$.
We may assume that $\mathcal{T}$ contains every lattice point
of $\partial P$ as a vertex (see, for example, \cite{AthVectors}). 
By Lemma \ref{upper}, if we set 
\begin{equation}\label{sunrise}
a(t) = \sum_{F \in \mathcal{T}} B_{F}(t)h_{F}(t) 
\end{equation}
\begin{equation}\label{sunset}
b(t) = t^{-l}\sum_{F \in \mathcal{T}} B_{F}'(t)h_{F}(t), 
\end{equation}
then $\bar{\delta}_{P}(t) = a(t) + t^{l}b(t)$. Since $mP$ contains no lattice points in its relative interior
for $m = 1, \ldots, l -1$, if $v$ lies in $\BOX(\sigma_{F}')\cap N'$ for some face $F$ of $\mathcal{T}$, then
$u(v) \geq l$. We conclude that $b(t)$ is a polynomial. By Lemma \ref{patriots}, the coefficients of $b(t)$ are non-negative
integers. 
Since every lattice point of $\partial P$ is a vertex of $\mathcal{T}$, 
if $v$ lies in $\BOX(\sigma_{F})\cap N'$ for some non-empty face $F$ of $\mathcal{T}$, then $u(v) \geq 2$. If we write 
$a(t) =  h_{\emptyset}(t) + t^{2}\sum_{F \in \mathcal{T}, F \neq \emptyset} t^{-2}B_{F}(t)h_{F}(t)$  then Lemma \ref{patriots}
implies that
$1 = a_{0} \leq a_{1} \leq a_{i}$ for $i = 2, \ldots, d - 1$. By Lemma \ref{badgers}, we are left with verifying that
$a(t) = t^{d}a(t^{-1})$ and $b(t) = t^{d - l}b(t^{-1})$. Using Lemmas \ref{patriots} and \ref{ball}, we compute
\[
t^{d}a(t^{-1}) = t^{d} \sum_{F \in \mathcal{T}} B_{F}(t^{-1})h_{F}(t^{-1}) 
= \sum_{F \in \mathcal{T}} B_{F}(t) t^{d - \dim F - 1} h_{F}(t^{-1}) 
= a(t),
\]
\[
t^{d - l}b(t^{-1}) = t^{d - l} t^{l} \sum_{F \in \mathcal{T}} B_{F}'(t^{-1})h_{F}(t^{-1})
= t^{-l}\sum_{F \in \mathcal{T}} B_{F}'(t) t^{d - \dim F - 1} h_{F}(t^{-1}) = b(t).
\]
\end{proof}

\begin{rem}
It follows from the above theorem that expressions (\ref{sunrise}) and (\ref{sunset}) are independent of the choice of lattice triangulation $\mathcal{T}$ and the choice of $\bar{v}$ in $lP \smallsetminus \partial(lP)$. 
\end{rem}

\begin{rem}
Let $K$ be the \emph{pyramid} over $\partial P$. That is, $K$ is the
truncation of $\partial \sigma$ at level $1$ and can be written as 
\[
K := \{ (x, \lambda) \in (N \times \mathbb{Z})_{\mathbb{R}} \mid x \in \partial(\lambda P),  0 < \lambda \leq 1  \} \cup \{ 0 \}.
\]
We may view $K$ as a polyhedral complex and consider its Ehrhart polynomial $f_{K}(m)$ and associated
Ehrhart $\delta$-polynomial $\delta_{K}(t)$ (p1 \cite{HibStar}, Chapter XI \cite{HibAlgebraic}). 
By the proof of Theorem \ref{goblue}, 
\[
a(t)/(1 - t)^{d + 1} = \sum_{F \in \mathcal{T}} \sum_{ \stackrel{v \in \sigma \cap N'}{\{ v \} \in \BOX(\sigma_{F})} } (1 + t + \cdots + t^{l-1})t^{u(v)}.
\]
It follows from Remark \ref{quo} that $a(t)/(1 - t)^{d + 1}$ is the generating series of $f_{K}(m)$ and hence that 
$a(t) = \delta_{K}(t)$.
With the terminology of \cite{HibStar}, $K$ is \emph{star-shaped} with respect to the origin and the fact that 
$1 = a_{0} \leq a_{1} \leq a_{i}$,
for $i = 2, \ldots, d - 1$, is a consequence of Hibi's results in \cite{HibStar}. 
\end{rem}


\begin{rem}\label{bonjovi}
By Theorem \ref{goblue}, $a_{0} \leq a_{1} \leq a_{2}$ and hence the coefficients of $a(t)$ are unimodal for $d \leq 5$ (c.f. Remark \ref{barnes}). 
\end{rem}

As a corollary, we obtain a combinatorial proof of a result of Stanley  \cite{StaHilbert2}. 
Recall, from Remark \ref{smh}, that a lattice polytope $P$ is reflexive 
if and only if it contains the origin in its 
relative interior
and, for every positive integer $m$, every non-zero lattice point in $mP$ lies on $\partial(n P)$ for a unique positive integer 
$n \leq m$. 

\begin{cor}\label{tata}
If $P$ is a lattice polytope of degree $s$ and codegree $l$, then $\delta_{P}(t) = t^{s} \delta_{P}(t^{-1})$ if and only if 
$lP$ is a translate of a reflexive polytope.  
\end{cor}
\begin{proof}
Since $t^{d}\bar{\delta}_{P}(t^{-1}) = t^{s}\delta_{P}(t^{-1}) (1 + t + \cdots + t^{l -1})$, we see that 
$\delta_{P}(t) = t^{s} \delta_{P}(t^{-1})$ if and only if $t^{d}\bar{\delta}_{P}(t^{-1}) = \bar{\delta}_{P}(t)$. By Lemma \ref{badgers},
we need to show that $b(t) = 0$ if and only if $lP$ is a translate of a reflexive polytope. 
By Remark \ref{quo}, 
$b(t) = 0$ if and only 
if every element $v$ in $\sigma \cap (N \times l\mathbb{Z})$
can be written as the sum of an element of $\partial(mlP) \times \{ ml \}$ and $m'(\bar{v}, l)$,
for some non-negative integers $m$ and $m'$. 
That is, $b(t)= 0$ if and only if $lP - \bar{v}$ is a reflexive polytope.  
\end{proof}

We now prove our second main result. 

\begin{theorem}\label{yoke}
Let $P$ be a $d$-dimensional lattice polytope of degree $s$ and codegree $l$.  The Ehrhart $\delta$-vector $(\delta_{0}, \ldots, \delta_{d})$
of $P$ satisfies the following inequalities.
\[
\delta_{1} \geq \delta_{d},
\]
\begin{equation*}
\delta_{2} + \cdots  + \delta_{i + 1} \geq \delta_{d - 1} + \cdots + \delta_{d - i}  \textrm{ for }     i = 0 , \ldots, \lfloor d/2     \rfloor - 1,
\end{equation*}
\begin{equation*}
\delta_{0} + \delta_{1} + \cdots + \delta_{i} \leq \delta_{s} + \delta_{s - 1} + \cdots + \delta_{s - i}
\textrm{ for }      i = 0, \ldots, d,
\end{equation*}
\begin{equation*}
\delta_{2 - l}  + \cdots + \delta_{0} +  \delta_{1} \leq \delta_{i} + \delta_{i - 1} + \cdots + \delta_{i - l + 1}  \textrm{ for }      
i = 2, \ldots, d - 1.
\end{equation*}
\end{theorem}
\begin{proof}
We observed in the introduction that $\delta_{1} \geq \delta_{d}$.  By Lemma \ref{magnify}, the second inequality is equivalent to
$a_{1} \leq a_{i}$, for $i = 2, \ldots, d - 1$, and 
the third inequality is equivalent to $b_{i} \geq 0$ for all $i$. Hence these inequalities follow from 
Theorem \ref{goblue}. When $l \geq 2$, the conditions above imply that $\bar{\delta}_{1} \leq \bar{\delta}_{i}$ for $i = 2, \ldots, d - 1$. By Lemma \ref{magnify}, this proves the final inequality when $l > 1$. When $l = 1$, the last inequality is Hibi's result (\ref{ineq3}).
\end{proof}

A lattice triangulation $\mathcal{T}$ of $\partial P$ is \emph{unimodular} if for every non-empty face $F$ of $\mathcal{T}$, the cone over $F \times \{1 \}$ in 
$N_{\mathbb{R}}'$
is non-singular. Equivalenty, $\mathcal{T}$ is unimodular if and only if $\BOX(\sigma_{F}) = \emptyset$, for every non-empty face $F$
of $\mathcal{T}$.

\begin{theorem}\label{dwarf}
Let $P$ be a $d$-dimensional lattice polytope. If $\partial P$ admits a regular unimodular lattice triangulation, then 
\[\delta_{i + 1} \geq \delta_{d - i} \] 
\[ \delta_{0} + \cdots + \delta_{i + 1} \leq \delta_{d} + \cdots + \delta_{d - i} + \binom{\delta_{1} - \delta_{d} + i + 1}{i + 1}, \]
for 
$i = 0, \ldots, \lfloor d/2 \rfloor - 1$.
\end{theorem}
\begin{proof}
{} From the above discussion, if $\partial P$ admits a regular unimodular triangulation then 
$B_{F}(t) = 0$ for every non-empty face $F$ of $\mathcal{T}$. By (\ref{sunrise}),
$a(t) = h_{F}(t)$, where $F$ is the empty face of 
$\mathcal{T}$. By Lemma \ref{badgers}, 
$1 = a_{0} \leq a_{1} \leq \cdots \leq a_{\lfloor d/2 \rfloor}$ and
$a_{i} \leq \binom{a_{1} + i - 1}{i}$
for $i = 1, \ldots, d$. The result now follows from the expression 
$a_{i + 1} = \delta_{0} + \cdots + \delta_{i + 1} - \delta_{d} - \cdots - \delta_{d - i}$ (see (\ref{coeff1})).
\end{proof}

\begin{rem}
When $P$ is the regular simplex of dimension $d$, $\delta_{P}(t) = 1$ and both inequalities in Theorem \ref{dwarf} are equalities.  
\end{rem}

\begin{rem}
Recall that if $P$ is a reflexive polytope then the coefficients of the Ehrhart $\delta$-vector are symmetric (Remark \ref{smh}). In this case, Theorem \ref{dwarf} implies that 
if $\partial P$ admits a regular, unimodular lattice triangulation then the coefficients of the $\delta$-vector are symmetric and unimodal. 
Note that if $P$ is reflexive then $P$ admits a regular unimodular lattice triangulation if and only if $\partial P$ admits a regular unimodular lattice triangulation.
Hence this result is a consequence of the theorem of Athanasiadis stated in the introduction (Theorem 1.3 \cite{AthVectors}).
This special case was first proved by Hibi in \cite{HibEhrhart}. 
\end{rem}

\begin{rem}
Recall that $\delta_{1} = |P \cap N| - (d + 1)$ and $\delta_{d} = |(P \smallsetminus \partial P) \cap N|$.
Hence $\delta_{1} = \delta_{d}$ if and only if $|\partial P \cap N| = d + 1$. If $|\partial P \cap N| = d + 1$ and $\partial P$ admits a regular, unimodular lattice triangulation
then, by Theorems  \ref{yoke} and \ref{dwarf}, $\delta_{0} + \cdots + \delta_{i + 1} = \delta_{d} + \cdots + \delta_{d - i} + 1$,
for 
$i = 0, \ldots, \lfloor d/2 \rfloor - 1$. This implies that $\delta_{i + 1} = \delta_{d - i}$ for 
$i = 0, \ldots, \lfloor d/2 \rfloor - 1$. If, in addition, $P$ is reflexive then the coefficients of the $\delta$-vector are symmetric and hence 
$\delta_{P}(t) = 1 + t + \cdots + t^{d}$.
\end{rem}

\begin{rem}\label{goop}
Let $P$ be a lattice polytope of dimension $d$ in $N$ and let $Q$ be the convex hull of $P \times \{1\}$ and the origin in $(N \times \mathbb{Z})_{\mathbb{R}}$. That is, $Q$ is the pyramid 
over $P$. If $\partial Q$ admits a regular unimodular lattice triangulation $\mathcal{T}$ then $\mathcal{T}$ restricts to give regular unimodular lattice triangulations of $P$ and $\partial P$.
\end{rem}

\begin{rem}
Hibi gave an example of a $4$-dimensional reflexive lattice polytope whose boundary does not admit a regular unimodular 
lattice triangulation  (Example 36.4 \cite{HibAlgebraic}). 
By Remark \ref{goop}, 
there are examples of $d$-dimensional lattice polytopes $P$ such that $\partial P$ does not admit a regular unimodular lattice triangulation for $d \geq 4$. On the other hand, if $P$ is a lattice polytope 
of dimension $d \leq 3$, then $\partial P$ always admits a regular, unimodular lattice triangulation. In fact, any regular triangulation of $\partial P$ containing every lattice point as a vertex is 
necessarily unimodular. This follows from the fact that if $Q$ is a lattice polytope of dimension $d' \leq 2$ and $|Q \cap N| = d' + 1$, then $Q$ is isomorphic to the regular $d'$-simplex. 
\end{rem}


\bibliographystyle{amsplain}
\bibliography{alan}

\end{document}